\documentclass[a4paper,11pt]{amsart}
\usepackage{amsmath,amssymb,amsthm}
\usepackage{fontspec}
\usepackage{unicode-math}
\usepackage{graphicx,enumitem}
\usepackage{booktabs,tabularx}
\usepackage[margin=27mm,headheight=15pt]{geometry}
\usepackage{placeins}
\usepackage{tikz}
\usetikzlibrary{arrows.meta,decorations.markings,decorations.pathreplacing}
\usepackage{hyperref}
\hypersetup{hidelinks,pdfauthor={Noboru Ito}}
\newcommand{\Q}{\mathbb Q}
\newcommand{\Z}{\mathbb Z}
\newcommand{\St}{\mathrm{St}}
\newcommand{\W}{\mathsf W}
\newcommand{\B}{\mathsf B}
\newcommand{\dd}{\mathrm d}
\DeclareMathOperator{\ind}{ind}
\DeclareMathOperator{\Tot}{Tot}
\DeclareMathOperator{\im}{Im}
\numberwithin{equation}{section}

\newtheorem{theorem}{Theorem}[section]
\newtheorem{reductionlemma}{Lemma}

\newtheorem{corollary}[theorem]{Corollary}
\theoremstyle{definition}
\newtheorem{fact}[theorem]{Fact}
\theoremstyle{remark}
\newtheorem{remark}[theorem]{Remark}
\title[Spectral sequences from \(\St\) to \(J\)]{Spectral sequences from \(\St\) homology\texorpdfstring{\\}{ }recovering \(J^-\) and \(J^+\)}
\author{Noboru Ito}
\date{29 November 2026}

\address{Department of Mathematics, Faculty of Engineering, Shinshu University,
Wakasato 4-17-1, Nagano, Nagano 380-8553, Japan}
\email{nito@shinshu-u.ac.jp}

\begin{document}

\begin{abstract}

Arnold's invariants are \(\St\), \(J^-\), and \(J^+=J^-+\text{the number of double points}\). The quantization polynomials of \(\St\) include the vertex state sum \(\St_q^v\) of \cite{Ito23} and the edge state sum \(\St_q\) of \cite{Ito26St}. The edge state sum is used for categorification. The quantization polynomial \(P_C(q)\) of \(J^-\) is due to Viro \cite{Viro96}. By the relation \(J^+=J^-+\text{the number of double points}\), the quantization polynomial of \(J^+\) is obtained, in essence, from \(P_C(q)\) together with the quantization polynomial \(V_q\) of the number of double points \cite{LP13}. These categorifications were obtained in \cite{Ito26St,Ito26}. We call the homological lifts of \(\St_q\), \(P_C(q)\), and \(V_q\) St homology, \(J^-\) homology, and double-point-count \(V\) homology, respectively. The pair of the latter two, which gives the homological lift of \(J^+\), is called \(J^+\) homology. The main result can then be stated as follows. Using de Rham and Koszul complexes, we construct two spectral sequences of families of complexes containing St homology, converging to \(J^-\) homology and to the image of the coefficient sequences obtained from double-point-count \(V\) homology, respectively, and obtain \(J^+\) homology. The spectral sequences thus directly relate the homological lifts of the quantization polynomials of Arnold's three invariants.

\end{abstract}\maketitle\enlargethispage{3pt}

\section{Introduction}

Arnold's invariants are \(\St\), \(J^-\), and \(J^+=J^-+\text{the number of double points}\).
For \(\St\), the vertex state sum is due to Ito \cite[Theorem~1 and Section~5]{Ito23}, and the edge state sum used for categorification is given in \cite[Definition~1.1]{Ito26St} (the notation \(\St_q^v,\St_q\) follows \cite[Theorem~(A)]{Ito26St}).
\[
 \St_q^v(C)=\sum_{v\in V(C)} w(v)q^{\ind(v)},\qquad
 \St_q(C)=\St_q^e(C)=\sum_{e\in E(C)} s(e)q^{\ind(e)}.
\]
\(C\) is a closed curve based on an exterior edge or the corresponding long curve \cite{Ito26St}, \(V(C)\) is the vertex set, and \(E(C)\) is the edge set of \(C\).
For \(s(e),w(v)\), the edge and vertex signs defined by intersection numbers, see \cite{Ito23,Ito26}. Below, \(\ind\) is the average at an edge or vertex of the Alexander indices of the adjacent faces \cite{Viro96}.

Consider the smoothing of the curve \(C\) with orientation \(o\). Let \(R_j(C,o)\) denote the faces of Alexander index \(j\) in its complement. The quantization polynomial \(P_C(q)\) of \(J^-\) is then given by \cite{Viro96} (for notation, cf. \cite[(1)]{Ito26}):
\[
 P_C(q)=\sum_{j\in\Z}\chi\bigl(R_j(C,o)\bigr)q^j.
\]

Together with \(P_C(q)\), the following polynomial
\[
 V_q(C)=\sum_{v\in V(C)} q^{\ind(v)}\bigl(q^{1/2}-q^{-1/2}\bigr)
\]
gives the \(J^+\) quantization polynomial \(I_q\) \cite{LP13}, as in \cite[Equations~(I3), (I4)]{Ito26}.\footnote{In \cite{Ito26}, subscripts indicating the orientation were used; here they are omitted.}
Indeed,
{\small
\[
 I_q(C)=I_q(\widetilde C)-\tfrac12 V_q(C),\quad
 P_C(q)=1+\bigl(q^{1/2}-q^{-1/2}\bigr)I_q(\widetilde C),\quad
 I_q(\widetilde C)=\sum_n\operatorname{rot}(\widetilde C_n)q^{\ind(\widetilde C_n)}.
\]}
Here \(\widetilde C_n\in\pi_0(\widetilde C)\) satisfy \(\widetilde C=\bigsqcup_n\widetilde C_n\), \(\operatorname{rot}(\widetilde C_n)\) is the rotation number \(\pm1\), and \(\ind\) is the average Alexander index of the adjacent faces.

In this paper, we call the homological lift of \(V_q\) \cite[Section~3.6]{Ito26} double-point-count \(V\) homology, referring to what it categorifies (or \(J\) edge homology, referring to its construction).

The categorification of \(P_C(q)\), the quantization polynomial of \(J^-\), has two main advantages \cite[Theorem~1.1(H) and Propositions~4.1 and 5.1]{Ito26}. It yields a strictly stronger means of distinguishing curves than the polynomial. Moreover, Viro introduced \(P_C(q)\) as an analogue of quantum knot invariants \cite[Section~5.3]{Viro96}, and decategorification of the homological lift naturally yields a face state-sum model with greater distinguishing power than \(P_C(q)\). The categorification of \(V_q\) gives a homological lift of \(J^+\) as the pair of circle homology and double-point-count \(V\) homology \cite[Theorem~1.1(E) and Section~3.6]{Ito26}. Decategorification also recovers the classical relation with \(P_C(q)\) \cite{LP13}, \cite[Equation~(I5)]{Ito26}.
The categorification of \(\St_q^{e}(C)\) (written simply as \(\St_q(C)\) in this paper unless otherwise specified) is likewise strictly stronger at distinguishing curves than the polynomial \cite[Theorem~(A) and Section~1.7]{Ito26St}. The construction of its categorification also leads naturally to the edge state sum as a quantization polynomial obtained by decategorification. The direct geometric relationship between \(\St\) and either \(J^-\) or \(J^+\) nevertheless remained unclear.

In this paper, we construct two spectral sequences of families of complexes containing St homology, converging to \(J^-\) homology and to the image of the coefficient sequences obtained from double-point-count \(V\) homology, respectively, and obtain \(J^+\) homology. The spectral sequences thus directly relate the homological lifts of the quantization polynomials of Arnold's three invariants.

The constructions are as follows.
We first construct a family of spectral sequences for which the direct sum of the zeroth pages contains St weighted circle homology, a homological lift of the quantization polynomial of Arnold's strangeness. With the appropriate coefficients, the direct sum of the terminal terms is identified with the circle homology, a homological lift of the quantization polynomial of Arnold's \(J^-\).
We next construct a spectral sequence for the quantization polynomial of the number of double points. When the reduced curve has crossings, its zeroth page contains the coefficient-sequence image obtained from St edge homology (via Lemma~\ref{lem:st-wb}), and the sequence converges to the image of the coefficient sequences obtained from double-point-count \(V\) homology (for \(2\)-reduction, see \hyperref[sec:terminal-supplement]{\emph{Supplement on the terminal term}} and Lemma~\ref{lem:reduction-count}). These constructions directly relate the quantization polynomials of Arnold's three invariants through spectral sequences.

Table~\ref{tableTerm} summarizes the terminology used in this paper and in the cited references.
The homologies in \cite[\S\S~3.1, 3.2, and 3.6]{Ito26} and \cite[\S1.6 and Remark~1.14]{Ito26St} have integer coefficients. In this paper, we construct spectral sequences in parallel in two settings, using the de Rham complex and the Koszul complex. The de Rham construction uses polynomial rings over \(\Q\), so in this setting we also use rational coefficients for the homologies in the preceding papers in place of integer coefficients. The Koszul construction requires no such change to the coefficients used in the preceding papers.

\begin{table}[!tbp]
\centering
\caption{Summary of the terminology used in this paper and in the cited references. The dash in the \emph{Notation here} column indicates that no notation is needed in this paper.}\label{tableTerm}
{\footnotesize
\setlength{\tabcolsep}{2pt}
\renewcommand{\arraystretch}{1.02}
\setlength{\aboverulesep}{1.5pt}
\setlength{\belowrulesep}{1.5pt}
\begin{tabularx}{\linewidth}{@{}
  >{\raggedright\arraybackslash}p{.117\linewidth}
  >{\raggedright\arraybackslash}X
  >{\raggedright\arraybackslash}p{.122\linewidth}
  >{\centering\arraybackslash}p{.125\linewidth}
  >{\raggedright\arraybackslash}p{.148\linewidth}
  >{\raggedright\arraybackslash}p{.144\linewidth}
  >{\raggedright\arraybackslash}p{.184\linewidth}
  @{}}
\toprule

\textbf{Homology} & \textbf{Type} & \textbf{De\-cate\-go\-rified}\newline\textbf{polynomial}
 & \textbf{Notation here} & \textbf{Notation in source} & \textbf{Definition in source}
 & \textbf{Spectral-sequence}\newline\textbf{term} \\
\midrule
St homology & St weighted circle homology & \(\St_q^e(=\St_q)\)
 & \(H^{i,j,k}(C_b^r)\) & \(H^{i,j,k}(C^r)\) & \cite{Ito26St}\newline \S1.6, after (1.4)
 & Each cycle included in a direct-sum component of \(E_0\) \\
 & St edge homology & \(\St_q^e(=\St_q)\)
 & \(H_{\St}\) & \(H_1(\mathcal E(C))\) & \cite{Ito26St}\newline Remark~1.14
 & \\
 & WB-labeled St edge homology & \(\St_q^{\mathrm{WB}}\)\newline \eqref{eq:wb-character}
 & \(H_{\mathrm{WB}}\) & \textemdash & Equation~\eqref{eq:edge-cycles}\newline Concept introduced in this paper
 & Coefficient-sequence image included in the direct sum of \(E_0\) \\
\midrule
\(J^-\) homology & circle homology & \(I_q(\widetilde C_o)\)
 & \(\mathcal H_{\mathrm{circle}}^{i,j,\varepsilon}(C,o)\) & \(\mathcal H_{\mathrm{circle}}^{i,j,\varepsilon}(C,o)\)
 & \cite{Ito26}, Thm.\newline 1.1(B), \S3.2
 & Identified with the direct sum of \(E_\infty\) \\
 & region homology & \(P_C(q)\)
 & \textemdash & \(\mathcal H_{\mathrm{region}}^{i,j}(C,o)\)
 & \cite{Ito26}, Thm.\newline 1.1(A), \S3.1 & \\
\midrule
double-point-count \(V\) homology & \(J\) edge homology & \(V_q(C)\)
 & \(H_J\) & \(\mathcal H_{\mathrm{edge}}(C,o)\) & \cite{Ito26}\newline \S3.6
 & \textsuperscript{1} \\
\midrule
\(J^+\) homology & \textsuperscript{2}
 & \(I_q(C)\) & \textemdash & \textemdash & Concept introduced in this paper & \\

\bottomrule
\end{tabularx}\par
\smallskip
\begin{tabularx}{\linewidth}{@{}l@{\hspace{.55em}}>{\raggedright\arraybackslash}X@{}}
\textsuperscript{1} & The image of the coefficient sequences of the WB canonical cycles is identified with the module spanned by the coefficient sequences read off at each crossing from the monochromatic words in \(E_\infty\). \\
\textsuperscript{2} & Pair of circle homology and double-point-count \(V\) homology
\end{tabularx}\par}
\end{table}

\section{From St homology to \texorpdfstring{\(J^-\)}{J-minus} homology}

\label{sec:circle}
\begin{theorem}\label{thm:circle}

Let \(C_b\) be a based closed curve in the sense of \cite[Section~1]{Ito26St}. Let \(C_b^r\) be its \(2\)-reduction and \(C^r\) the curve obtained by forgetting the base point of \(C_b^r\).
Use the traversal numbering from \(-\infty\) on the long curve corresponding to \(C_b^r\) (obtained by an opening that decreases the rotation number by \(1\)). When the number of crossings is \(m\geq 1\), the edges are numbered
\begin{equation}\label{eq:edge-order}
 e_1,\ldots,e_{2m},e_{2m+1}=e_\infty .
\end{equation}
The edge set used for St weighted circle homology is \(\{e_1,\ldots,e_{2m}\}\). When \(m\geq1\), order the circles of the smoothing of \(C^r\) by the least edge number on each circle. Denote the configuration by
\[
 (\widetilde{C^r}_1,\ldots,\widetilde{C^r}_n).
\]
There is a family
\(\{E(\widetilde{C^r}_1),\ldots,E(\widetilde{C^r}_n)\}\) of spectral sequences arising from the de Rham and Koszul complexes. The family has the following properties.

\begin{enumerate}[label=(\arabic*),leftmargin=*]

\item For the spectral sequences arising from the de Rham complexes, the direct sum of the zeroth pages contains the following \(\Q\)-linear direct sum of the cycles of St weighted circle homology.

\[
 \Q\langle[z_1]\rangle\oplus\cdots\oplus\Q\langle[z_n]\rangle .
\]
Each direct summand \(\Q\langle[z_i]\rangle\) corresponds to the respective \(\widetilde{C^r}_i\).

For the spectral sequences arising from the Koszul complexes, the direct sum of the zeroth pages contains the corresponding direct sum over \(\Z\).

\[
 \Z\langle[z_1]\rangle\oplus\cdots\oplus\Z\langle[z_n]\rangle .
\]

\item Let \(N\) be the number of edges of the \(i\)th circle \(\widetilde{C^r}_i\). Let \(K=K_i\) be the de Rham complex for \(\Q[x_1,\ldots,x_N]\). The terminal term of the spectral sequence arising from the de Rham complex is isomorphic to \(H^0(K)\).
Let \(K=K_i\) be the Koszul complex of the review in Section~\ref{sec:review}. The terminal term of the spectral sequence arising from the Koszul complex is isomorphic to \(H^N(K)\).
This does not apply when the coefficients are changed so that the de Rham and Koszul complexes are isomorphic.
Also, if \(C_b^r\) has no crossings, take \(N=1\) for its circle.

\item The all-\(\W\) classes of the terminal terms of the spectral sequences arising from the de Rham complexes completely give the circle homology obtained by smoothing the curve \(C^r\) along its orientation. On each \(\widetilde{C^r}_i\), the terminal term is \(E_\infty=\Q\langle\{\W\text{-word}\}\rangle\).
The all-\(\B\) classes of the terminal terms of the spectral sequences arising from the Koszul complexes completely give the circle homology obtained by smoothing along \((C^r,\bar o)\), which reverses only the orientation of \(C^r=(C^r,o)\). On each \(\widetilde{C^r}_i\), the terminal term is \(E_\infty=\Z\langle\{\B\text{-word}\}\rangle\).
It is immediate that nontriviality in degree \(i=1\) and \(j=\text{Alexander index}\) of circle homology are recovered.

\end{enumerate}\end{theorem}

\subsection*{Step 1. Outline}
Consider the de Rham complex with the two-colour words of length \(N\), in the letters \(\B,\W\), as a basis over the polynomial ring. In the word description of the differential, one \(\W\) changes to \(\B\). The horizontal differential \(d_h\) changes odd positions, while the vertical differential \(d_v\) changes even positions. The terminal term is generated by a single monochromatic word. In the construction below, we present the spectral sequence arising from the de Rham complex first and give the Koszul case alongside it where the constructions proceed in parallel.

The details are as follows.

\begin{enumerate}[label=(\arabic*),leftmargin=*]

\item For the \(i\)th circle, let \(N\) be its number of edges. The complex \(K=K_i\) is the de Rham complex for \(\Q[x_1,\ldots,x_N]\), or the Koszul complex of the review in Section~\ref{sec:review}. Use the column filtration of the total complex. The bigraded modules of \(K\) form the zeroth page \(E_0^{p,q}=K^{p,q}\).

When \(C_b^r\) has no crossings, take \(N=1\). The single letter \(\W\) corresponds to agreement with the orientation of the circle, and \(\B\) to the opposite orientation. In this case \(J^+(C^r)=J^-(C^r)\).

\item The horizontal and vertical differentials are

\[
 d_h:K^{p,q}\longrightarrow K^{p+1,q},
 \qquad d_v:K^{p,q}\longrightarrow K^{p,q+1}.
\]

The differential on the zeroth page is \(d_0=d_v\).
\item The first page is

\[
 E_1^{p,q}=H^q(K^{p,*},d_v).
\]

\item The second page satisfies

\[
 E_2^{p,q}=H^p(E_1^{*,q},d_1)\cong E_\infty^{p,q}.
\]

Here \(d_1\) is induced by \(d_h\).

\end{enumerate}

\subsection*{Step 2. The St-state and the WB-state}
The state in \cite{Ito26St} records the configuration of Seifert circles with the vertex positions. It induces signs on the vertices and edges. St weighted circle homology uses the state of the \(2\)-reduced curve. In Steps~2--4 we abbreviate \(C_b^r\) or \(C^r\) as \(C\) when there is no ambiguity.

For each smoothed circle \(\widetilde C_i\), denote its generator by \([\widetilde C_i]\) and consider the free \(\Q\)-module

\begin{equation}\label{eq:circle-module}
 \bigoplus_i\Q\langle[\widetilde C_i]\rangle .
\end{equation}

Use the edge numbering in \eqref{eq:edge-order}. On a circle with numbered edges there is a least edge number, which is used to order the \(\widetilde C_i\).

\begin{figure}[htbp]
\centering
\begingroup
\input{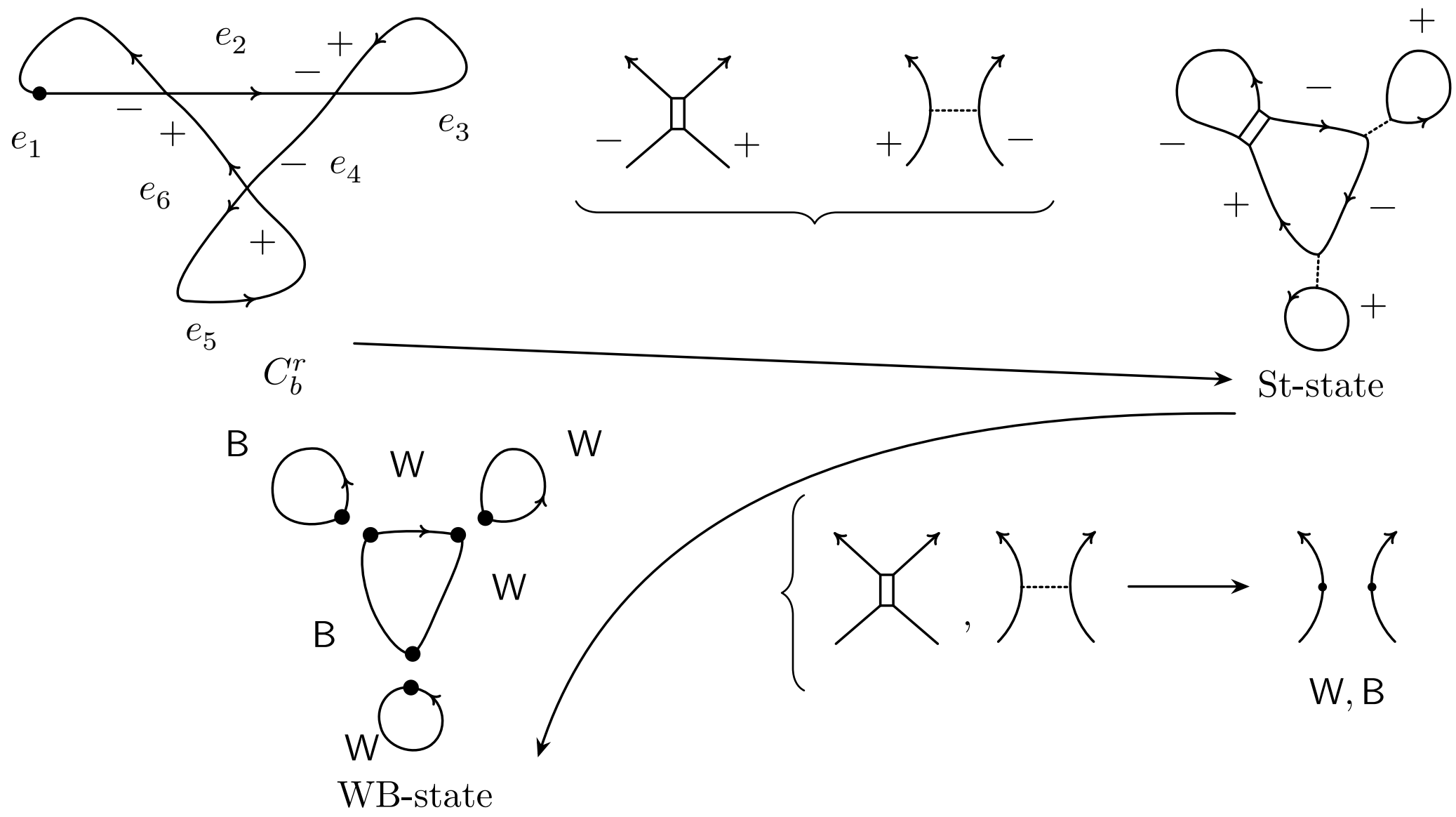}

\newcommand{\FromStToJWideEdgeDetail}{%
  \draw (-.48,-.48) -- (-.061,-.12) -- (-.061,.17)
        (.48,-.48) -- (.061,-.12) -- (.061,.17);
  \draw (-.061,-.12) -- (.061,-.12) (-.061,.17) -- (.061,.17);
  \draw[-{To[length=3.3pt,width=3.8pt]}] (-.061,.17) -- (-.50,.57);
  \draw[-{To[length=3.3pt,width=3.8pt]}] (.061,.17) -- (.50,.57);
}
\newcommand{\FromStToJSmoothedDetail}{%
  \draw[-{To[length=3.3pt,width=3.8pt]}] (-.45,-.46)
     .. controls (-.18,-.19) and (-.13,.29) .. (-.47,.58);
  \draw[-{To[length=3.3pt,width=3.8pt]}] (.45,-.46)
     .. controls (.18,-.19) and (.13,.29) .. (.47,.58);
}
\newcommand{\FromStToJSingleEdgeDetail}{%
  \FromStToJSmoothedDetail
  \draw[densely dotted] (-.233,.055) -- (.233,.055);
}

\begin{tikzpicture}[
  x=1cm, y=1cm, font=\small, fromsttoj contour,
  every node/.style={inner sep=1pt}
]
  \begin{scope}[shift={(2.24,5.84)},scale=1.46]
    \FromStToJBasedCurve{postaction={decorate},decoration={markings,
      mark=at position .15 with {\arrow{To[length=3.2pt,width=3.8pt]}},
      mark=at position .38 with {\arrow{To[length=3.2pt,width=3.8pt]}},
      mark=at position .53 with {\arrow{To[length=3.2pt,width=3.8pt]}},
      mark=at position .65 with {\arrow{To[length=3.2pt,width=3.8pt]}},
      mark=at position .78 with {\arrow{To[length=3.2pt,width=3.8pt]}},
      mark=at position .88 with {\arrow{To[length=3.2pt,width=3.8pt]}}}}
    \node at (-1.42,.29) {\(e_1\)};
    \node at (-.10,.96) {\(e_2\)};
    \node at (1.34,.40) {\(e_3\)};
    \node at (.64,.13) {\(e_4\)};
    \node at (-.29,-.95) {\(e_5\)};
    \node at (-.59,-.04) {\(e_6\)};
    \node at (-.76,.49) {\(-\)};
    \node at (.39,.73) {\(-\)};
    \node at (.60,.94) {\(+\)};
    \node at (.30,.135) {\(-\)};
    \node at (.10,-.34) {\(+\)};
    \node at (-.48,.36) {\(+\)};
  \end{scope}
  \node at (2.60,4.08) {\(C_b^r\)};

  \begin{scope}[shift={(6.30,6.53)}]
    \FromStToJWideEdgeDetail
    \node at (-.65,-.27) {\(-\)};
    \node at (.65,-.27) {\(+\)};
  \end{scope}
  \begin{scope}[shift={(8.91,6.53)}]
    \FromStToJSingleEdgeDetail
    \node at (-.62,-.26) {\(+\)};
    \node at (.62,-.26) {\(-\)};
  \end{scope}
  \draw[line width=.5pt,decorate,decoration={brace,amplitude=6pt}]
    (9.84,5.73) -- (5.34,5.73);
  \draw[fromsttoj flow] (3.25,4.39) -- (11.53,4.05);

  \begin{scope}[shift={(12.36,5.91)},scale=1.77]
    \FromStToJStState
    \node at (-.79,.18) {\(-\)};
    \node at (-.01,.48) {\(-\)};
    \node at (.54,.86) {\(+\)};
    \node at (.33,-.16) {\(-\)};
    \node at (.28,-.66) {\(+\)};
    \node at (-.45,-.12) {\(+\)};
  \end{scope}
  \node at (12.36,4.01) {St-state};

  \begin{scope}[shift={(3.82,2.15)},scale=1.78]
    \FromStToJCircleConfiguration{fromsttoj orientation=.17}{fromsttoj orientation=.165}
                                 {fromsttoj orientation=.32}{fromsttoj orientation=.90}
    \FromStToJMarkedPoints
    \node at (-.94,.73) {\(\B\)};
    \node at (-.04,.62) {\(\W\)};
    \node at (.90,.73) {\(\W\)};
    \node at (.50,-.03) {\(\W\)};
    \node at (-.28,-.88) {\(\W\)};
    \node at (-.48,-.28) {\(\B\)};
  \end{scope}
  \node at (3.82,.14) {WB-state};

  \draw[fromsttoj flow] (11.55,3.73)
    .. controls (8.50,3.76) and (5.83,3.05) .. (4.98,.49);
  \draw[line width=.5pt,decorate,decoration={brace,amplitude=6pt}]
    (7.49,1.12) -- (7.49,2.96);
  \begin{scope}[shift={(8.27,2.04)}]
    \FromStToJWideEdgeDetail
  \end{scope}
  \node at (9.03,1.72) {,};
  \begin{scope}[shift={(9.77,2.04)}]
    \FromStToJSingleEdgeDetail
  \end{scope}
  \draw[fromsttoj flow] (10.55,2.10) -- (11.69,2.10);
  \begin{scope}[shift={(12.61,2.04)}]
    \FromStToJSmoothedDetail
    \fill (-.233,.055) circle[radius=1.18pt]
          (.233,.055) circle[radius=1.18pt];
  \end{scope}
  \node at (12.61,1.08) {\(\W,\B\)};
\end{tikzpicture}
\endgroup
\caption{A curve \(C_b^r\) and its St- and WB-states. The base point is omitted in the St- and WB-states.}
\label{fig:wb}
\end{figure}

Since we start from the state of \cite{Ito26St}, the data in \([\widetilde C]\) include the circle configuration, the circle orientations, vertex positions, edge numbers, the base point of the original curve, and the edge signs. In this paper, the base-point marking is omitted where no ambiguity arises (Figures~\ref{fig:wb} and~\ref{fig:example-flow}). We call this state the St-state and write the module in \eqref{eq:circle-module} in terms of its cycles as

\begin{equation*}
 \bigoplus_i\Q\langle[\widetilde C_i]_{\St}\rangle .
\end{equation*}

For each \(\widetilde C_i\), compare the sign of an edge on it with the rotation sign of the circle, taking \(+1\) for counterclockwise orientation and \(-1\) for clockwise orientation. Assign \(\W\) when the signs agree and \(\B\) when they differ. On a circle with numbered edges, the letters read in the direction of the circle's orientation form a WB word whose first letter corresponds to the least numbered edge.

A state obtained by assigning the letters \(\W,\B\) to the edges is called a WB-state. We assume that both the St- and WB-states carry the orientation of each Seifert circle. Thus, for a given curve \(C\), the St-state and WB-state correspond one-to-one, and in terms of the WB-state we write \eqref{eq:circle-module} as

\begin{equation}\label{eq:circle-module-wb}
 \bigoplus_i\Q\langle[\widetilde C_i]_{\mathrm{WB}}\rangle .
\end{equation}

See Figure~\ref{fig:wb}. In Step~2, the Koszul case uses the same construction with \(\Z\) in place of \(\Q\).

\subsection*{Step 3. The complex for each circle}
For each circle in \eqref{eq:circle-module-wb}, let \(N\) be its number of edges, using the convention of Step~1 when there are no crossings. Take all \(2^N\) words of length \(N\) in \(\W,\B\). In the de Rham case, interpret \(\W\) as \(1\) and \(\B\) in the \(i\)th position as \(\dd x_i\), and set

\begin{equation*}
 A=\Q[x_1,\ldots,x_N],\qquad C^k=\Omega^k_{A/\Q}.
\end{equation*}

In the Koszul case, using \(\#\B(w)=\#\{i:w_i=\B\}\), set

\begin{equation*}
 R=\Z[t_1^{\pm1},\ldots,t_N^{\pm1}],\quad f_i=t_i-1,\quad
 K_R^k=\bigoplus_{\#\B(w)=k}R\,w.
\end{equation*}

The spectral sequences of these complexes are outlined in Step~1 and described in detail in Section~\ref{sec:review}. Choosing the all-\(\W\) word on every circle corresponds to \([\widetilde C]\). Choosing the all-\(\B\) word on every circle corresponds to \([(\widetilde C,\bar o)]\), obtained from \((C,\bar o)\), where \(\bar o\) reverses the orientation \(o\) of \(C\).

\subsection*{Step 4. The terminal term}
On each \(\widetilde C_i\), the terminal term is
\begin{equation*}
 E_\infty=\Q\langle\{\W\text{-word}\}\rangle .
\end{equation*}
This completely gives the circle homology obtained by smoothing the curve \(C=C^r\) along its orientation.

For the spectral sequence of the Koszul complex, the terminal term on each \(\widetilde C_i\) is
\begin{equation*}
 E_\infty=\Z\langle\{\B\text{-word}\}\rangle .
\end{equation*}
This completely gives the circle homology obtained by smoothing along \((C^r,\bar o)\), which reverses only the orientation of the curve \(C=(C^r,o)\).
It is immediate that nontriviality in degree \(i=1\) and \(j=\text{Alexander index}\) of circle homology are recovered.
Thus the zeroth pages contain the cycle lines of the St weighted circle homology of \cite{Ito26St}, and the terminal terms give the cycle lines of the circle homology of \cite{Ito26}, with the circle configuration.

\subsection*{Supplement on the terminal term}\label{sec:terminal-supplement}
In both the de Rham and Koszul cases, the circle homology at the terminal term is that of the \(2\)-reduction \(C^r\) of the given curve \(C\). The notation \(C^r\) is from \cite{Ito26St}. To recover \(J^-(C)\) itself by decategorification, the correction is \(2n\), where \(n\) is the minimum number of strong \(2\)-moves required for the \(2\)-reduction.

The \(2\)-reduction repeatedly uses crossing-decreasing \(2\)-moves until there are no \(2\)-gons to remove. Its result is unique by \cite[Proposition~1.4]{Ito26St}. For related reduction arguments, see \cite[Theorem~2.2 (2)]{IT13}, \cite[Section~2]{Khov97}, and \cite[Section~5, Notation~3 and Fact~2]{Ito26R}. The difference \(J^-(C^r)-J^-(C)=|J^-(C^r)-J^-(C)|\) is twice the number of crossing-decreasing strong \(2\)-moves, denoted \(s2b\), in the reduction. Such a move increases \(J^-\) by \(2\).

\begin{reductionlemma}\label{lem:reduction-count}

The minimum number \(n\) of strong \(2\)-moves required to obtain the \(2\)-reduction \(C^r\) of \(C\) is uniquely determined.

\end{reductionlemma}\begin{proof}

Let \(h\) be a \(2\)-reduction from \(C\) to \(C^r\), and let \(n(h,s2b)\) be the number of \(s2b\) moves in \(h\).
Under Arnold's defining rules for \(J^-\), an \(s2b\) move increases \(J^-\) by \(2\), and weak \(2\)-moves and \(3\)-moves leave it unchanged \cite{Arn94}. Regular homotopies are described by \(2\)- and \(3\)-moves. Consequently, along the reduction \(h\),

\begin{equation*}
 J^-(C^r)-J^-(C)=2n(h,s2b).
\end{equation*}

For two \(2\)-reductions \(h'\) and \(h''\), this gives

\[
 n(h',s2b)=\frac{J^-(C^r)-J^-(C)}{2}=n(h'',s2b).
\]

This proves the assertion.
\end{proof}

\section{From St homology to double-point-count \texorpdfstring{\(V\)}{V} homology}

\label{sec:edges}

This section concerns double-point-count \(V\) homology, the categorification of the quantization polynomial \(V_q(C^r)=2\bigl(I_q(\widetilde{C^r})-I_q(C^r)\bigr)\) of the number of double points \(J^+-J^-\). First suppose that \(C^r\) has \(N>0\) crossings. Use the based curve \(C_b^r\) and the numbering in \eqref{eq:edge-order}, with \(m=N\). Recall the signs in St edge homology from \cite{Ito26St}: when the intersection number of \(v\) is negative, \(s(e_l^v)=-1,s(e_r^v)=+1\); when it is positive, \(s(e_l^v)=+1,s(e_r^v)=-1\).

Use the edge complex for St from \cite[Remark~1.14]{Ito26St} and the edge complex for \(J^+\) from \cite[Section~3.6]{Ito26}, replacing the integer coefficients by a polynomial ring over \(\Q\) for the de Rham case and by a Laurent polynomial ring over \(\Z\) for the Koszul case.
The canonical cycle in St edge homology is \(\bigl[s(e_l^v)(e_l^v-e_r^v)\bigr]\). For an edge \(e\) on the Seifert circle \(\widetilde{C^r}_i\), let \(\operatorname{rot}(\widetilde{C^r}_i)\in\{\pm1\}\) be its rotation number. Define a new sign \(t(e)\) by
\begin{equation}\label{eq:wb-edge-sign}
 \begin{cases}
 \text{if }s(e)=\operatorname{rot}(\widetilde{C^r}_i),\text{ label }\W\text{ and set }t(e)=+1,\\
 \text{if }s(e)\ne\operatorname{rot}(\widetilde{C^r}_i),\text{ label }\B\text{ and set }t(e)=-1.
 \end{cases}
\end{equation}
Preserving the isomorphism of St edge homology \(H_{\St}\), there is a one-to-one correspondence
\begin{equation}\label{eq:edge-cycles}
 \bigl[s(e_l^v)(e_l^v-e_r^v)\bigr]
 \longleftrightarrow
 \bigl[t(e_l^v)(e_l^v-e_r^v)\bigr].
\end{equation}
Call this WB-labeled St edge homology \(H_{\mathrm{WB}}\), and call \(\bigl[t(e_l^v)(e_l^v-e_r^v)\bigr]\) the WB canonical cycle.

\begin{reductionlemma}\label{lem:st-wb}
For a fixed curve \(C=C_b^r\),
\[
 \bigl[s(e_l^v)(e_l^v-e_r^v)\bigr]
 \longmapsto
 \bigl[t(e_l^v)(e_l^v-e_r^v)\bigr]
\]
induces an isomorphism \(H_{\St}(C)\cong H_{\mathrm{WB}}(C)\).
\end{reductionlemma}

\begin{proof}
By the definition in \eqref{eq:wb-edge-sign}, the signs \(s(e)\) in the St-state and \(t(e)\) in the WB-state correspond one-to-one.
\end{proof}

Write \(ch_q^{\St}\) for the polynomial obtained by decategorifying using the canonical cycles \(\bigl[s(e_l^v)(e_l^v-e_r^v)\bigr]\) defined by \(s(e)\):
\begin{equation}\label{eq:st-character}
 \sum_v s(e_l^v)\bigl(q^{\ind(e_l^v)}-q^{\ind(e_r^v)}\bigr).
\end{equation}
The canonical cycle in the edge homology \(H_J\) of \(J^+\) homology (in particular, the homology whose decategorification contributes the correction term of \(I_q\)) is \(\bigl[e_l^v-e_r^v\bigr]\).

Consider a WB-state in which the WB labels are chosen freely, independently of the signs \(s(e)\) fixed by \(C=C_b^r\), and define
\begin{equation}\label{eq:wb-character}
 ch_q^{\mathrm{WB}}(C)=\sum_v t(e_l^v)\bigl(q^{\ind(e_l^v)}-q^{\ind(e_r^v)}\bigr).
\end{equation}
Recall that in Section~\ref{sec:circle} we assigned WB labels using the signs \(s(e)\) of the curve. Here we assign \(\W\) and \(\B\) independently of \(s(e)\). If we set \(t(v)=-1\) (\(s(v)=-1\), resp.) when \(e_l^v\) has the label \(\B\) (\(s(e_l^v)=-1\), resp.), and \(t(v)=1\) (\(s(v)=+1\), resp.) otherwise, then for the vertex set \(\{v_1,\ldots,v_N\}\), the respective coefficient sequences of the edge complexes corresponding to \(ch_q^{\St}\) and to \(ch_q^{\mathrm{WB}}\) are
\begin{equation}\label{eq:coefficient-sequences}
 (s(v_1),\ldots,s(v_N)),\qquad (t(v_1),\ldots,t(v_N)).
\end{equation}
Lemma~\ref{lem:st-wb} describes the correspondence between the two coefficient sequences.

Their characters give \(\St_q^e\) on the \(E_0\) page and the contribution of the number of double points to \(I_q\) on the \(E_\infty\) page, respectively. Indeed, the former is clear from the defining equation \eqref{eq:st-character}, and in the latter case, the sum in the correction term of \cite[Equation~(I3)]{Ito26} can be written, using the indices in \cite[Section~3.6]{Ito26}, as follows:

\begin{equation*}
 V_q(C)=\sum_v\bigl(q^{\ind(e_l^v)}-q^{\ind(e_r^v)}\bigr).
\end{equation*}

\begin{corollary}\label{cor:wb-vq}
\[
 t(e_l^v)=1\;(\forall v)
 \quad\Longrightarrow\quad
 ch_q^{\mathrm{WB}}(C)=V_q(C).
\]
\end{corollary}

The sequence of signs \(t(v_i)\) obtained from the first coefficient sequence of signs \(s(v_i)\) in \eqref{eq:coefficient-sequences} by the correspondence of Lemma~\ref{lem:st-wb} gives words in \(\W,\B\). The coefficient sequence at \(E_\infty\) corresponds to the all-\(\W\) word in the de Rham case and to the all-\(\B\) word in the Koszul case. Take all words of length \(N_{i}\) for the edge coefficient sequences on each circle \(\widetilde{C^r}_i\), and form the spectral sequences of their de Rham and Koszul complexes as in Section~\ref{sec:review}. Let us be more explicit.

In the de Rham case, the word corresponding to the signed coefficient sequence gives the following embedding into the zeroth page.

\begin{equation*}
 \im(H_{\mathrm{WB}}\longrightarrow\Q^N)
 =\Q\langle(t(v_1),\ldots,t(v_N))\rangle
 \hookrightarrow E_0 .
\end{equation*}

At the terminal term, the all-\(\W\) class corresponds to the coefficient sequence \((1,\ldots,1)\) from double-point-count \(V\) homology.

\begin{equation*}
 E_\infty\cong\Q\langle[\W^N]\rangle
 \cong\im(H_J\longrightarrow\Q^N)
 =\Q\langle(1,\ldots,1)\rangle .
\end{equation*}

In the Koszul case, the corresponding embedding into the zeroth page is

\begin{equation*}
 \im(H_{\mathrm{WB}}\longrightarrow\Z^N)
 =\Z\langle(t(v_1),\ldots,t(v_N))\rangle
 \hookrightarrow E_0 .
\end{equation*}

At the terminal term, the all-\(\B\) class corresponds to the coefficient sequence \((-1,\ldots,-1)\) from double-point-count \(V\) homology for \((C^r,\bar o)\).

\begin{equation*}
 E_\infty\cong\Z\langle[\B^N]\rangle
 \cong\im(H_J(C^r,\bar o)\longrightarrow\Z^N)
 =\Z\langle(-1,\ldots,-1)\rangle .
\end{equation*}

For \(N>0\), these spectral sequences therefore converge to the respective coefficient-sequence images obtained from double-point-count \(V\) homology.

\noindent When \(N=0\), for the de Rham complex, the zeroth and terminal pages in bidegree \((0,0)\) are \(\Q\), so we define both St edge homology and double-point-count \(V\) homology to be \(\Q\) (note that these are taken to be empty in \cite{Ito26,Ito26St}). The character \(ch_q^{\mathrm{WB}}\) is zero by \eqref{eq:wb-character}.
For the Koszul complex, the zeroth and terminal pages in bidegree \((0,0)\) are \(\Z\), so we define both St edge homology and double-point-count \(V\) homology to be \(\Z\) (note that these are taken to be empty in \cite{Ito26,Ito26St}). The character \(ch_q^{\mathrm{WB}}\) is zero by \eqref{eq:wb-character}.

In summary, as in Theorem~\ref{thm:circle}, we have the following.

\begin{theorem}\label{thm:double-point-count}
Let \(C_b\) be a based closed curve in the sense of \cite[Section~1]{Ito26St}. Let \(C_b^r\) be its \(2\)-reduction, and let \(C^r\) be the curve obtained by forgetting the base point of \(C_b^r\).
Number the edges by following the long curve corresponding to \(C_b^r\) (obtained by an opening that decreases the rotation number by \(1\)) from \(-\infty\) in the direction of its orientation. When the number of crossings is \(N\geq1\), the edges are numbered
\[
 e_1,\ldots,e_{2N},\ e_{2N+1}=e_\infty .
\]
The edge set used is \(\{e_1,\ldots,e_{2N}\}\). When \(N\geq1\), order the smoothing circles of \(C^r\) by the smallest edge number on each circle. Write the circle configuration as
\[
 (\widetilde{C^r}_1,\ldots,\widetilde{C^r}_n).
\]
There are families of spectral sequences \(\{E(\widetilde{C^r}_1),\ldots,E(\widetilde{C^r}_n)\}\), arising from the de Rham and Koszul complexes, with the following properties.

\begin{enumerate}[label=(\arabic*),leftmargin=*]
\item The direct sum of the zeroth pages of the spectral sequences arising from the de Rham complexes contains the following direct sum over \(\Q\), obtained from the cycles of St weighted circle homology:
\[
 \Q\langle[z_1]\rangle\oplus\cdots\oplus\Q\langle[z_n]\rangle .
\]
For each \(i\), the summand \(\Q\langle[z_i]\rangle\) corresponds to \(\widetilde{C^r}_i\).
The direct sum of the zeroth pages of the spectral sequences arising from the Koszul complexes contains the corresponding direct sum over \(\Z\):
\[
 \Z\langle[z_1]\rangle\oplus\cdots\oplus\Z\langle[z_n]\rangle .
\]

\item Let \(N_i\) be the number of edges of the \(i\)th circle \(\widetilde{C^r}_i\). Let \(K=K_i\) be the de Rham complex of \(\Q[x_1,\ldots,x_{N_i}]\). The terminal term of the spectral sequence arising from the de Rham complex is isomorphic to \(H^0(K)\).
Let \(K=K_i\) be the Koszul complex of the review in Section~\ref{sec:review}. The terminal term of the spectral sequence arising from the Koszul complex is isomorphic to \(H^{N_i}(K)\).
This does not apply when the coefficients are changed as in Remark~\ref{rem:coefficient-change} so that the de Rham and Koszul complexes are isomorphic.

There is also a construction compatible with decategorification when the number of crossings is zero.

\item The all-\(\W\) classes of the terminal terms of the spectral sequences arising from the de Rham complexes give double-point-count \(V\) homology. The decategorification to the quantization polynomial of the number of double points is given by \(ch_q^{\mathrm{WB}}\) in \eqref{eq:wb-character}.
On each \(\widetilde{C^r}_i\), the terminal term is \(E_\infty=\Q\langle\{\W\text{-word}\}\rangle\).
The all-\(\B\) classes of the terminal terms of the spectral sequences arising from the Koszul complexes give double-point-count \(V\) homology. In this case the decategorification to the quantization polynomial of the number of double points is given by \(-ch_q^{\mathrm{WB}}\).
On each \(\widetilde{C^r}_i\), the terminal term is \(E_\infty=\Z\langle\{\B\text{-word}\}\rangle\).
\end{enumerate}
\end{theorem}

\section{An example}\label{sec:example}
In this section, we calculate St weighted circle homology for the specific curve in Figure~\ref{fig:wb}, calculate the spectral sequences starting with the corresponding family of \(E_0\) pages, and obtain \(E_\infty\). We then confirm that the image of the coefficient sequence of the resulting \(E_\infty\) terms is (as one of the \(J^-\) homologies) circle homology.

We also examine an example of the spectral sequence starting with a family of zeroth pages containing St edge homology (via Lemma~\ref{lem:st-wb}), whose module spanned by the coefficient sequences read off at each crossing from the monochromatic words in \(E_\infty\) is double-point-count \(V\) homology.
\begin{figure}[!htbp]
\centering
\newcommand{\FromStToJStateReference}{(Figure~\ref{fig:wb})}
\begingroup
\input{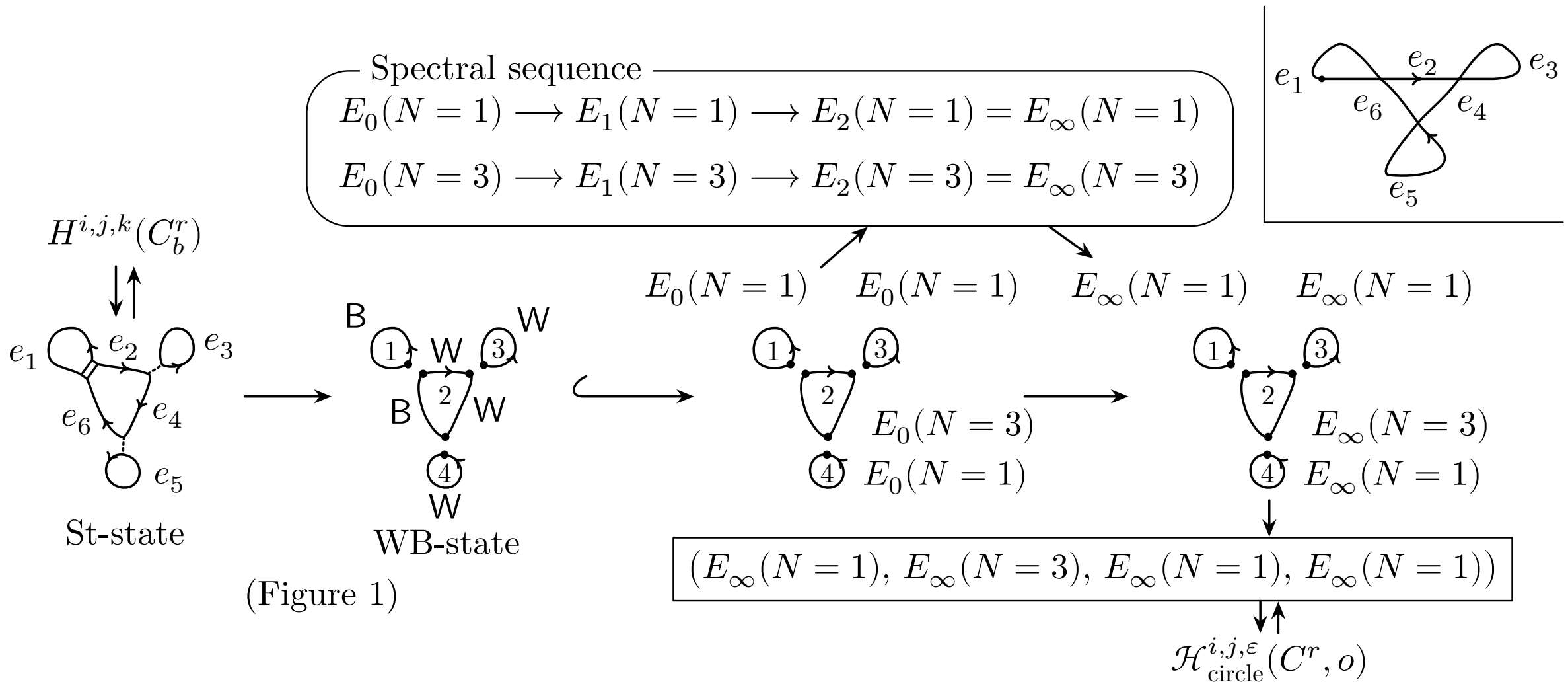}
\newcommand{\FromStToJCircleNumbers}{%
  \node[circle number] at (-.555,.473) {1};
  \node[circle number] at (-.012,.030) {2};
  \node[circle number] at (.521,.497) {3};
  \node[circle number] at (-.018,-.749) {4};
}
\newcommand{\FromStToJPageLabels}[1]{%
  \node[anchor=south east] at (-.15,.88) {\(#1(N=1)\)};
  \node[anchor=south west] at (.21,.88) {\(#1(N=1)\)};
  \node[anchor=west] at (.37,-.28) {\(#1(N=3)\)};
  \node[anchor=west] at (.29,-.78) {\(#1(N=1)\)};
}
\begin{tikzpicture}[
  x=1cm, y=1cm, font=\small, fromsttoj contour,
  every node/.style={inner sep=1pt},
  circle number/.style={font=\scriptsize, inner sep=0pt}
]
  \begin{scope}[shift={(14.14,5.23)},scale=.72]
    \FromStToJBasedCurve{postaction={decorate},decoration={markings,
      mark=at position .146 with {\arrow{To[length=3pt,width=3.8pt]}},
      mark=at position .741 with {\arrow{To[length=3pt,width=3.8pt]}}}}
    \node at (-1.78,.59) {\(e_1\)};
    \node at (.07,.84) {\(e_2\)};
    \node at (1.76,.80) {\(e_3\)};
    \node at (.76,.22) {\(e_4\)};
    \node at (-.19,-.98) {\(e_5\)};
    \node at (-.69,.22) {\(e_6\)};
  \end{scope}
  \draw[line width=.5pt] (12.60,6.40) -- (12.60,4.24) -- (15.59,4.24);

  \node[draw,rounded corners=12pt,line width=.45pt,inner xsep=9pt,inner ysep=7pt]
    (page-flow) at (7.65,4.98) {%
      \(\begin{gathered}
       E_0(N=1)\longrightarrow E_1(N=1)\longrightarrow E_2(N=1)=E_\infty(N=1)\\[3pt]
       E_0(N=3)\longrightarrow E_1(N=3)\longrightarrow E_2(N=3)=E_\infty(N=3)
      \end{gathered}\)%
    };
  \node[anchor=west,fill=white,inner xsep=3pt,inner ysep=1pt]
    at ([xshift=15pt]page-flow.north west) {Spectral sequence};
  \draw[fromsttoj flow] (8.16,3.79) -- (8.63,0 |- page-flow.south);
  \draw[fromsttoj flow] (10.45,0 |- page-flow.south) -- (10.89,3.88);

  \node at (1.19,4.12) {\(H^{i,j,k}(C_b^r)\)};
  \draw[fromsttoj flow] (1.10,3.80) -- (1.10,3.29);
  \draw[fromsttoj flow] (1.28,3.29) -- (1.28,3.80);

  \begin{scope}[shift={(1.19,2.48)}]
    \FromStToJStState
    \node[anchor=east] at (-.80,.41) {\(e_1\)};
    \node at (-.01,.51) {\(e_2\)};
    \node[anchor=west] at (.74,.52) {\(e_3\)};
    \node at (.43,-.17) {\(e_4\)};
    \node at (-.48,-.24) {\(e_6\)};
    \node[anchor=west] at (.24,-.82) {\(e_5\)};
    \node at (0,-1.34) {St-state};
  \end{scope}

  \draw[fromsttoj flow] (2.39,2.50) -- (3.26,2.50);

  \begin{scope}[shift={(4.41,2.48)}]
    \FromStToJCircleConfiguration{fromsttoj orientation=.17}{fromsttoj orientation=.165}
                                 {fromsttoj orientation=.32}{fromsttoj orientation=.90}
    \FromStToJMarkedPoints
    \FromStToJCircleNumbers
    \node at (-.91,.83) {\(\B\)};
    \node at (.00,.48) {\(\W\)};
    \node at (.87,.79) {\(\W\)};
    \node at (.43,-.11) {\(\W\)};
    \node at (-.46,-.13) {\(\B\)};
    \node at (-.015,-1.08) {\(\W\)};
    \node at (0,-1.44) {WB-state};
  \end{scope}

  \node at (3.16,.50) {\FromStToJStateReference};

  \draw[fromsttoj flow] (5.81,2.70)
    .. controls (5.67,2.64) and (5.63,2.50) .. (5.78,2.50)
    -- (6.89,2.50);

  \begin{scope}[shift={(8.24,2.48)}]
    \FromStToJCircleConfiguration{fromsttoj orientation=.17}{fromsttoj orientation=.165}
                                 {fromsttoj orientation=.32}{fromsttoj orientation=.90}
    \FromStToJMarkedPoints
    \FromStToJCircleNumbers
    \FromStToJPageLabels{E_0}
  \end{scope}

  \draw[fromsttoj flow] (10.20,2.50) -- (11.24,2.50);

  \begin{scope}[shift={(12.65,2.48)}]
    \FromStToJCircleConfiguration{fromsttoj orientation=.17}{fromsttoj orientation=.165}
                                 {fromsttoj orientation=.32}{fromsttoj orientation=.90}
    \FromStToJMarkedPoints
    \FromStToJCircleNumbers
    \FromStToJPageLabels{E_\infty}
  \end{scope}

  \node[anchor=east,draw,line width=.45pt,inner xsep=4pt,inner ysep=3pt]
    (limit-list) at (15.09,.77)
    {\(\bigl(E_\infty(N=1),\, E_\infty(N=3),\, E_\infty(N=1),\, E_\infty(N=1)\bigr)\)};
  \draw[fromsttoj flow] (12.65,1.46) -- (12.65,0 |- limit-list.north);
  \draw[fromsttoj flow] (12.56,0 |- limit-list.south) -- (12.56,.13);
  \draw[fromsttoj flow] (12.74,.13) -- (12.74,0 |- limit-list.south);
  \node at (12.65,-.13) {\(\mathcal H_{\mathrm{circle}}^{i,j,\varepsilon}(C^r,o)\)};
\end{tikzpicture}
\endgroup
\caption{The flow from a given \(C_b^r\) through its St- and WB-states and the family of zeroth pages containing St weighted circle homology, and from the family of terminal terms to the coefficient sequences corresponding to circle homology. The base point is omitted in the St- and WB-states. The figure omits the all-\(\B\) words in the terminal terms of the Koszul case.}
\label{fig:example-flow}
\end{figure}

\subsection*{Example calculation (1)}
\(S=\Q\) (de Rham case) or \(\Z\) (Koszul case).
The words of length \(N\) corresponding to the six edges of \(C_b^r\) in Figure~\ref{fig:example-flow} are associated with the Seifert circles \(\widetilde C_1,\widetilde C_2,\widetilde C_3,\widetilde C_4\), respectively.
\[
\begin{aligned}
 \text{St weighted circle homology}\cong{}&
S\langle[\widetilde C_1]\rangle\quad\text{with grades }(i,j,k)=(1,\tfrac12,-1) \\
\oplus{}&S\langle[\widetilde C_2]\rangle\quad\text{with grades }(i,j,k)=(1,-\tfrac12,-1-1+1) \\
\oplus{}&S\langle[\widetilde C_3]\rangle\quad\text{with grades }(i,j,k)=(1,\tfrac12,1) \\
\oplus{}&S\langle[\widetilde C_4]\rangle\quad\text{with grades }(i,j,k)=(1,\tfrac12,1).
\end{aligned}
\]
\noindent\([\widetilde C_1]\) (the word \(\B\) from \(e_1\) labelled \(\B\)): \(N=1\).
\[
\begin{gathered}
S\langle[\widetilde C_1]\rangle\hookrightarrow E_0(N=1)\longrightarrow E_\infty(N=1):\\
\text{circle homology cycle }\langle[\widetilde C_1]\rangle\quad\text{with grades }(i,j,\varepsilon)\\
=\begin{cases}
  \text{de Rham case:} & (1,\tfrac12,1),\\
  \text{Koszul case:} & (1,-\tfrac12,1)
\end{cases}.
\end{gathered}
\]
\noindent\([\widetilde C_2]\) (the word \(\W\W\B\) from \(e_2\) labelled \(\W\), \(e_4\) labelled \(\W\), \(e_6\) labelled \(\B\)): \(N=3\).
\[
\begin{gathered}
S\langle[\widetilde C_2]\rangle\hookrightarrow E_0(N=3)\longrightarrow E_\infty(N=3):\\
\text{circle homology cycle }\langle[\widetilde C_2]\rangle\quad\text{with grades }(i,j,\varepsilon)\\
=\begin{cases}
  \text{de Rham case:} & (1,-\tfrac12,0),\\
  \text{Koszul case:} & (1,\tfrac12,1)
\end{cases}.
\end{gathered}
\]
\noindent\([\widetilde C_3]\) (the word \(\W\) from \(e_3\) labelled \(\W\)): \(N=1\).
\[
\begin{gathered}
S\langle[\widetilde C_3]\rangle\hookrightarrow E_0(N=1)\longrightarrow E_\infty(N=1):\\
\text{circle homology cycle }\langle[\widetilde C_3]\rangle\quad\text{with grades }(i,j,\varepsilon)\\
=\begin{cases}
  \text{de Rham case:} & (1,\tfrac12,1),\\
  \text{Koszul case:} & (1,-\tfrac12,1)
\end{cases}.
\end{gathered}
\]
\noindent\([\widetilde C_4]\) (the word \(\W\) from \(e_5\) labelled \(\W\)): \(N=1\).
\[
\begin{gathered}
S\langle[\widetilde C_4]\rangle\hookrightarrow E_0(N=1)\longrightarrow E_\infty(N=1):\\
\text{circle homology cycle }\langle[\widetilde C_4]\rangle\quad\text{with grades }(i,j,\varepsilon)\\
=\begin{cases}
  \text{de Rham case:} & (1,\tfrac12,1),\\
  \text{Koszul case:} & (1,-\tfrac12,1)
\end{cases}.
\end{gathered}
\]
Thus,
\par\noindent de Rham case
\[
\text{circle homology}\text{ for }(C^r,o)
\cong\Q\langle[\widetilde C_1]\rangle\oplus\Q\langle[\widetilde C_2]\rangle\oplus\Q\langle[\widetilde C_3]\rangle\oplus\Q\langle[\widetilde C_4]\rangle.
\]
\noindent Koszul case
\[
\text{circle homology}\text{ for }(C^r,\bar o)
\cong\Z\langle[\widetilde C_1]\rangle\oplus\Z\langle[\widetilde C_2]\rangle\oplus\Z\langle[\widetilde C_3]\rangle\oplus\Z\langle[\widetilde C_4]\rangle.
\]
These are the correspondences.
\subsection*{Example calculation (2)}
The cycles \([z]=\bigl[t(e_l^v)(e_l^v-e_r^v)\bigr]\) at the vertices \(v\) for WB-labeled St edge homology \(H_{\mathrm{WB}}\) are as follows.
\[
\begin{aligned}
v_1:\quad&\bigl(e_1\text{ with grades }(\B,\tfrac12),\ e_6\text{ with grades }(\B,-\tfrac12)\bigr),\\
&[z_1]=\bigl[(-1)(e_1-e_6)\bigr],\\
&ch_q^{\mathrm{WB}}([z_1])=(-1)\bigl(q^{1/2}-q^{-1/2}\bigr).
\end{aligned}
\]
\[
\begin{aligned}
v_2:\quad&\bigl(e_2\text{ with grades }(\W,-\tfrac12),\ e_3\text{ with grades }(\W,\tfrac12)\bigr),\\
&[z_2]=\bigl[(+1)(e_3-e_2)\bigr],\\
&ch_q^{\mathrm{WB}}([z_2])=(+1)\bigl(q^{1/2}-q^{-1/2}\bigr).
\end{aligned}
\]
\[
\begin{aligned}
v_3:\quad&\bigl(e_4\text{ with grades }(\W,-\tfrac12),\ e_5\text{ with grades }(\W,\tfrac12)\bigr),\\
&[z_3]=\bigl[(+1)(e_5-e_4)\bigr],\\
&ch_q^{\mathrm{WB}}([z_3])=(+1)\bigl(q^{1/2}-q^{-1/2}\bigr).
\end{aligned}
\]
The signs \((-1,1,1)\) determining the coefficients of WB-labeled St edge homology \(H_{\mathrm{WB}}\) correspond to \((\B,\W,\W)\).
The convergence calculation depends on the chosen complex.

\noindent de Rham case:\par\nopagebreak\noindent
The cycles \(\bigl[t(e_l^v)(e_l^v-e_r^v)\bigr]\) at the vertices \(v_i\) for double-point-count \(V\) homology are as follows.
\[
\begin{aligned}
v_1:\quad&\bigl(e_1\text{ with grades }(\W,\tfrac12),\ e_6\text{ with grades }(\W,-\tfrac12)\bigr),\\
&[z_1]=\bigl[(+1)(e_1-e_6)\bigr],\\
&ch_q^{\mathrm{WB}}([z_1])=(+1)\bigl(q^{1/2}-q^{-1/2}\bigr).
\end{aligned}
\]
\[
\begin{aligned}
v_2:\quad&\bigl(e_2\text{ with grades }(\W,-\tfrac12),\ e_3\text{ with grades }(\W,\tfrac12)\bigr),\\
&[z_2]=\bigl[(+1)(e_3-e_2)\bigr],\\
&ch_q^{\mathrm{WB}}([z_2])=(+1)\bigl(q^{1/2}-q^{-1/2}\bigr).
\end{aligned}
\]
\[
\begin{aligned}
v_3:\quad&\bigl(e_4\text{ with grades }(\W,-\tfrac12),\ e_5\text{ with grades }(\W,\tfrac12)\bigr),\\
&[z_3]=\bigl[(+1)(e_5-e_4)\bigr],\\
&ch_q^{\mathrm{WB}}([z_3])=(+1)\bigl(q^{1/2}-q^{-1/2}\bigr).
\end{aligned}
\]
This gives the cycles.
The signs \((1,1,1)\) determining the coefficients of WB-labeled \(J\) edge homology \(H_J\) correspond to \((\W,\W,\W)\). All \(t(e)\) are equal to \(1\).

\noindent Koszul case:\par\nopagebreak\noindent
The cycles \(\bigl[t(e_l^v)(e_l^v-e_r^v)\bigr]\) at the vertices \(v_i\) for double-point-count \(V\) homology are as follows.
\[
\begin{aligned}
v_1:\quad&\bigl(e_6\text{ with grades }(\B,\tfrac12),\ e_1\text{ with grades }(\B,-\tfrac12)\bigr),\\
&[z_1]=\bigl[(-1)(e_6-e_1)\bigr],\\
&ch_q^{\mathrm{WB}}([z_1])=(-1)\bigl(q^{1/2}-q^{-1/2}\bigr).
\end{aligned}
\]
\[
\begin{aligned}
v_2:\quad&\bigl(e_2\text{ with grades }(\B,\tfrac12),\ e_3\text{ with grades }(\B,-\tfrac12)\bigr),\\
&[z_2]=\bigl[(-1)(e_2-e_3)\bigr],\\
&ch_q^{\mathrm{WB}}([z_2])=(-1)\bigl(q^{1/2}-q^{-1/2}\bigr).
\end{aligned}
\]
\[
\begin{aligned}
v_3:\quad&\bigl(e_4\text{ with grades }(\B,\tfrac12),\ e_5\text{ with grades }(\B,-\tfrac12)\bigr),\\
&[z_3]=\bigl[(-1)(e_4-e_5)\bigr],\\
&ch_q^{\mathrm{WB}}([z_3])=(-1)\bigl(q^{1/2}-q^{-1/2}\bigr).
\end{aligned}
\]
The signs \((-1,-1,-1)\) determining the coefficients of WB-labeled \(J\) edge homology \(H_J\) correspond to \((\B,\B,\B)\). All \(t(e)\) are equal to \(-1\).
\par

\begin{table}[!htbp]
\centering
\caption{The spectral sequence \(E(N=1,3)\), de Rham complex case}
\label{tab:de-rham-example}
\begin{tabular}{@{}ccc@{}}
\toprule
 & \(E_1\) & \(E_2=E_\infty\) \\
\midrule
\(N=1\) & \(\begin{array}{c|cc} q\backslash p & 0 & 1\\\hline 0 & [\W] & [\B]\\ \end{array}\) & \(\begin{array}{c|cc} q\backslash p & 0 & 1\\\hline 0 & [\W] & 0\\ \end{array}\) \\[4pt]
\(N=3\) & \(\begin{array}{c|ccc} q\backslash p & 0 & 1 & 2\\\hline 0 & [\W\W\W] & [\B\W\W],\, [\W\W\B] & [\B\W\B]\\ 1 & 0 & 0 & 0\\ \end{array}\) & \(\begin{array}{c|ccc} q\backslash p & 0 & 1 & 2\\\hline 0 & [\W\W\W] & 0 & 0\\ 1 & 0 & 0 & 0\\ \end{array}\) \\[4pt]
\bottomrule
\end{tabular}
\end{table}

\begin{table}[!htbp]
\centering
\caption{The spectral sequence \(E(N=1,3)\), Koszul complex case}
\label{tab:koszul-example}
\begin{tabular}{@{}ccc@{}}
\toprule
 & \(E_1\) & \(E_2=E_\infty\) \\
\midrule
\(N=1\) & \(\begin{array}{c|cc} q\backslash p & 0 & 1\\\hline 0 & [\W] & [\B]\\ \end{array}\) & \(\begin{array}{c|cc} q\backslash p & 0 & 1\\\hline 0 & 0 & [\B]\\ \end{array}\) \\[4pt]
\(N=3\) & \(\begin{array}{c|ccc} q\backslash p & 0 & 1 & 2\\\hline 0 & 0 & 0 & 0\\ 1 & [\W\B\W] & [\B\B\W],\, [\W\B\B] & [\B\B\B]\\ \end{array}\) & \(\begin{array}{c|ccc} q\backslash p & 0 & 1 & 2\\\hline 0 & 0 & 0 & 0\\ 1 & 0 & 0 & [\B\B\B]\\ \end{array}\) \\[4pt]
\bottomrule
\end{tabular}
\end{table}

\FloatBarrier

\begin{figure}[!htbp]
\centering
\begingroup
\input{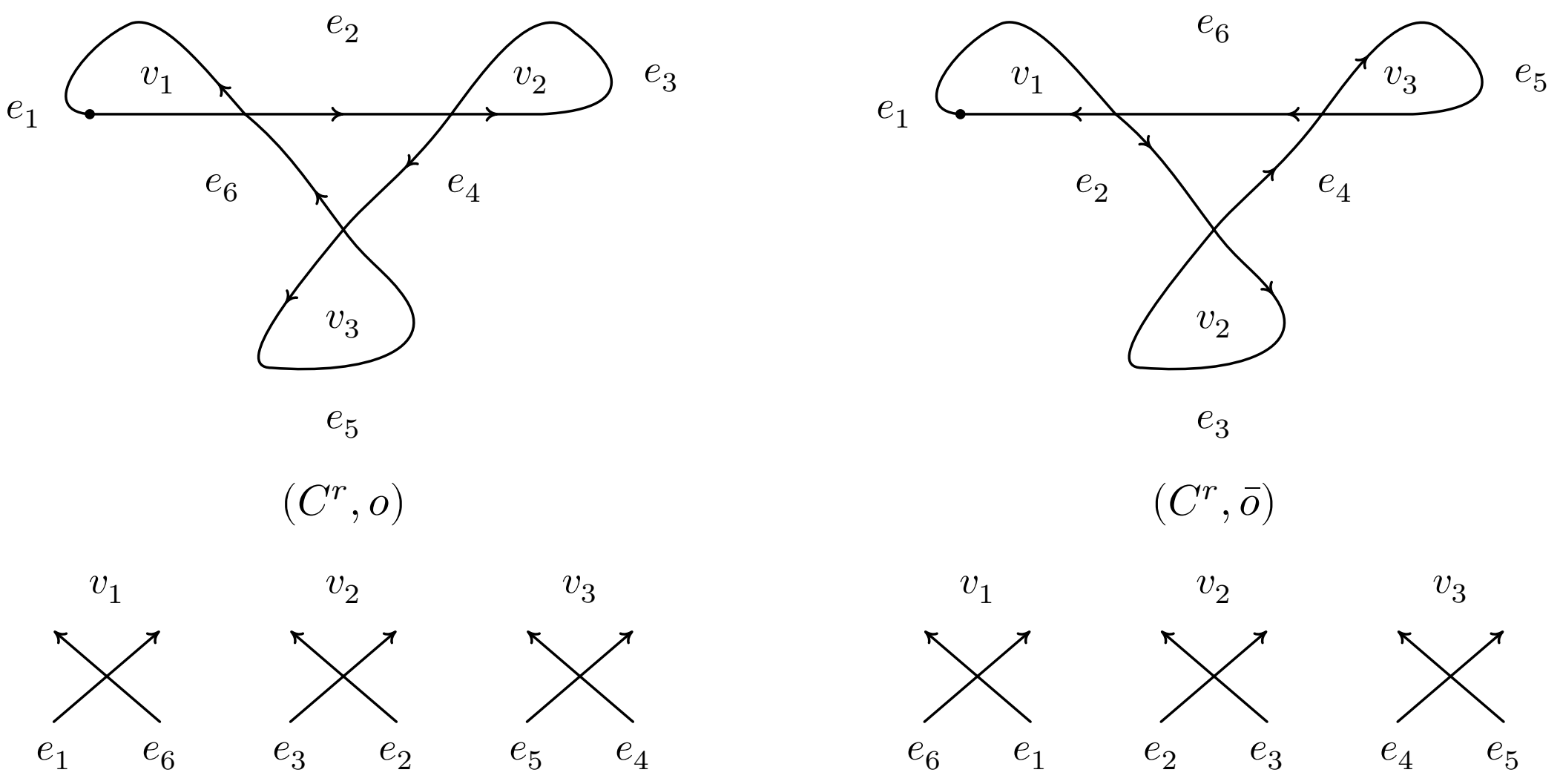}
\providecommand{\FromStToJFigureThreeRevisionColor}{black}
\newcommand{\FromStToJFigureThreeArrow}[1]{%
  #1{To[length=3.2pt,width=3.8pt,color=\FromStToJFigureThreeRevisionColor]}%
}
\newcommand{\FromStToJFigureThreeCrossing}[4]{%
  \begin{scope}[shift={(#1,-2.35)},color=\FromStToJFigureThreeRevisionColor]
    \node at (0,.46) {\(v_{#2}\)};
    \draw[-{To[length=3.2pt,width=3.8pt]}] (-.28,-.24) -- (.28,.24);
    \draw[-{To[length=3.2pt,width=3.8pt]}] (.28,-.24) -- (-.28,.24);
    \node at (-.28,-.42) {\(e_{#3}\)};
    \node at (.28,-.42) {\(e_{#4}\)};
  \end{scope}%
}
\resizebox{.90\linewidth}{!}{%
\begin{tikzpicture}[
  x=1.7cm, y=1.7cm, fromsttoj contour, font=\small,
  every node/.style={inner sep=1pt}
]
  \path[use as bounding box] (-1.86,-2.95) rectangle (6.46,1.23);

  \begin{scope}
    \FromStToJBasedCurve{postaction={decorate},decoration={markings,
      mark=at position .14 with {\arrow{To[length=3.2pt,width=3.8pt]}},
      mark=at position .225 with {\FromStToJFigureThreeArrow{\arrow}},
      mark=at position .46 with {\FromStToJFigureThreeArrow{\arrow}},
      mark=at position .56 with {\arrow{To[length=3.2pt,width=3.8pt]}},
      mark=at position .785 with {\FromStToJFigureThreeArrow{\arrow}},
      mark=at position .865 with {\FromStToJFigureThreeArrow{\arrow}}}}
    \node[text=\FromStToJFigureThreeRevisionColor] at (-1.69,.62) {\(e_1\)};
    \node at (0,1.07) {\(e_2\)};
    \node at (1.68,.81) {\(e_3\)};
    \node at (.64,.23) {\(e_4\)};
    \node at (0,-1.02) {\(e_5\)};
    \node at (-.64,.23) {\(e_6\)};
    \node at (-.98,.80) {\(v_1\)};
    \node at (.99,.80) {\(v_2\)};
    \node at (0,-.49) {\(v_3\)};
    \node[font=\normalsize] at (0,-1.44) {\((C^r,o)\)};

    \FromStToJFigureThreeCrossing{-1.25}{1}{1}{6}
    \FromStToJFigureThreeCrossing{0}{2}{3}{2}
    \FromStToJFigureThreeCrossing{1.25}{3}{5}{4}
  \end{scope}

  \begin{scope}[shift={(4.60,0)}]
    \FromStToJBasedCurve{postaction={decorate},decoration={markings,
      mark=at position .06 with {\FromStToJFigureThreeArrow{\arrowreversed}},
      mark=at position .18 with {\FromStToJFigureThreeArrow{\arrowreversed}},
      mark=at position .38 with {\arrowreversed{To[length=3.2pt,width=3.8pt]}},
      mark=at position .46 with {\FromStToJFigureThreeArrow{\arrowreversed}},
      mark=at position .71 with {\FromStToJFigureThreeArrow{\arrowreversed}},
      mark=at position .815 with {\arrowreversed{To[length=3.2pt,width=3.8pt]}}}}
    \node at (-1.69,.62) {\(e_1\)};
    \node at (-.64,.23) {\(e_2\)};
    \node at (0,-1.02) {\(e_3\)};
    \node at (.64,.23) {\(e_4\)};
    \node at (1.68,.81) {\(e_5\)};
    \node at (0,1.07) {\(e_6\)};
    \node at (-.98,.80) {\(v_1\)};
    \node at (0,-.49) {\(v_2\)};
    \node at (.99,.80) {\(v_3\)};
    \node[font=\normalsize] at (0,-1.44) {\((C^r,\bar{o})\)};

    \FromStToJFigureThreeCrossing{-1.25}{1}{6}{1}
    \FromStToJFigureThreeCrossing{0}{2}{2}{3}
    \FromStToJFigureThreeCrossing{1.25}{3}{4}{5}
  \end{scope}
\end{tikzpicture}%
}
\endgroup
\caption{\((C^r,o)\) and \((C^r,\bar{o})\) in Section~\ref{sec:example}(2).}
\label{fig:curve-orientations}
\end{figure}

\FloatBarrier

\section{Review of the de Rham and Koszul complexes}
We recall the basic facts about the de Rham and Koszul complexes used in this paper.

\label{sec:review}

\subsection{The de Rham complex}
Let \(A=\Q[x_1,\ldots,x_N]\). A word \(w=(w_1,\ldots,w_N)\) with \(w_i\in\{\W,\B\}\) specifies the positions \(i_1<\cdots<i_k\) at which the letter is \(\B\). This determines a bijection \(\iota\) from words to the exterior basis elements and its inverse:

\begin{equation*}
 \iota:w\longmapsto\dd x_{i_1}\wedge\cdots\wedge\dd x_{i_k},
 \qquad
 C^k:=\bigoplus_{\#\B(w)=k}A\,w\cong\Omega^k_{A/\Q}.
\end{equation*}

Let \(b_i(w)=\#\{j<i:w_j=\B\}\). The differential is

\begin{equation*}
 d(fw)=\sum_{i:w_i=\W}(-1)^{b_i(w)}
 \frac{\partial f}{\partial x_i}\,w^{(i)} .
\end{equation*}

Here \(w^{(i)}\) is obtained by changing the \(i\)th letter of \(w\) from \(\W\) to \(\B\). The complex \((C^\bullet,d)\) is the de Rham complex. Let \(\delta_i\) be the corresponding component of the differential, with value zero when the \(i\)th letter is \(\B\). Put

\begin{align*}
 d_h&=\sum_{i\ \mathrm{odd}}\delta_i,&
 d_v&=\sum_{i\ \mathrm{even}}\delta_i,\\
 p&=\#\B_{\mathrm{odd}}(w),&
 q&=\#\B_{\mathrm{even}}(w),&
 k&=p+q. 
\end{align*}
\begin{equation*}
 C^{p,q}=\bigoplus_{\substack{
     w:\ \#\B_{\mathrm{odd}}(w)=p\\
     \#\B_{\mathrm{even}}(w)=q
 }} A\,w .
\end{equation*}

The maps then have the bidegrees

\[
 d_h:C^{p,q}\longrightarrow C^{p+1,q},
 \qquad d_v:C^{p,q}\longrightarrow C^{p,q+1}.
\]

The total complex and column filtration are

\begin{equation*}
 \Tot^n C=\bigoplus_{p+q=n}C^{p,q},\qquad
 F^p\Tot^n C=\bigoplus_{\substack{i+j=n\\i\geq p}}C^{i,j}.
\end{equation*}

The differential satisfies \(d(F^p)\subseteq F^p\). Taking the vertical differential first, we have the pages

\begin{align*}
 E_0^{p,q}&=C^{p,q},\qquad d_0=d_v,\\
 E_1^{p,q}&=H^q(C^{p,*},d_v),\\
 E_2^{p,q}&=H^p(E_1^{*,q},d_1).
\end{align*}

Here \(d_1\) is induced by \(d_h\). The filtration is finite, so the limit is

\begin{equation*}
 E_\infty^{p,q}\cong
 F^pH^{p+q}(\Tot C)/F^{p+1}H^{p+q}(\Tot C).
\end{equation*}

For the standard double-complex spectral sequence and convergence, see \cite[Lemmas~12.25.1 and 12.25.3]{Stacks}. Let \(a=\lceil N/2\rceil\). By the algebraic Poincar\'e lemma \cite[Section~2, Lemma~2.1 and Remark~2.2]{Ant16},

\begin{equation*}
 E_1^{p,q}\cong
 \begin{cases}
   \Q[x_1,x_3,\ldots]^{\binom ap} & q=0,\\
   0 & q\ne0 .
 \end{cases}
\end{equation*}

There are \(\binom ap\) remaining words, each having \(\W\) at every even position. Applying the Poincar\'e lemma to the odd variables gives

\begin{equation*}
 E_2^{p,q}=E_\infty^{p,q}\cong
 \begin{cases}
   \Q\langle[\W^N]\rangle & (p,q)=(0,0),\\
   0 & \text{otherwise}.
 \end{cases}
\end{equation*}

\begin{corollary}\label{cor:de-rham}
The cohomology of the de Rham complex is

\begin{equation*}
 H^k(C^\bullet)\cong
 \begin{cases}
   \Q & k=0,\\
   0 & k\ne0.
 \end{cases}
\end{equation*}
\end{corollary}

\subsection{The Koszul complex}
Use the ring and word modules

\begin{equation*}
 R=\Z[t_1^{\pm1},\ldots,t_N^{\pm1}],\qquad
 f_i=t_i-1,\qquad
 K_R^k=\bigoplus_{\#\B(w)=k}R\,w .
\end{equation*}

with differential

\begin{equation*}
 d_K(rw)=\sum_{i:w_i=\W}(-1)^{b_i(w)}(t_i-1)r\,w^{(i)}
 \qquad (r\in R).
\end{equation*}

Let \(\delta_i\) be the component changing the \(i\)th letter from \(\W\) to \(\B\), with value zero when the \(i\)th letter is \(\B\). Use \(d_h=\sum_{i\ \mathrm{odd}}\delta_i\) and \(d_v=\sum_{i\ \mathrm{even}}\delta_i\). As in the de Rham case, \(p\) and \(q\) count the \(\B\)'s in odd and even positions, respectively. Set

\begin{align*}
 K_R^{p,q}
   &=\bigoplus_{\substack{
       w:\ \#\B_{\mathrm{odd}}(w)=p\\
       \#\B_{\mathrm{even}}(w)=q
      }} R\,w,\\
 F^p\Tot^n K_R
   &=\bigoplus_{\substack{i+j=n\\i\geq p}}K_R^{i,j}.
\end{align*}

The Koszul cohomology calculation for the regular sequence \(t_1-1,\ldots,t_N-1\) gives the following fact \cite[p.~75 and Corollary~4.8.4(i)]{Schapira}.
\begin{fact}\label{fact:koszul-h}

\begin{equation*}
 H^k(K_R^\bullet)\cong
 \begin{cases}
   R/(t_1-1,\ldots,t_N-1) & k=N,\\
   0 & k\ne N.
 \end{cases}
\end{equation*}
\end{fact}

\noindent\begin{minipage}{\linewidth}
\begin{fact}\label{fact:koszul-pages}
Let \(a=\lceil N/2\rceil\) and \(b=\lfloor N/2\rfloor\), and put

\begin{equation*}
 R_{\mathrm{odd}}
 =R/(t_2-1,t_4-1,\ldots)
 \cong\Z[t_1^{\pm1},t_3^{\pm1},\ldots].
\end{equation*}

For the spectral sequence taking the vertical differential first,

\begin{equation}\label{eq:koszul-e1}
 E_1^{p,q}\cong
 \begin{cases}
   R_{\mathrm{odd}}^{\binom ap} & q=b,\\
   0 & q\ne b,
 \end{cases}
\end{equation}
\begin{equation}\label{eq:koszul-e2}
 E_2^{p,q}=E_\infty^{p,q}\cong
 \begin{cases}
   \bigl(R/(t_1-1,\ldots,t_N-1)\bigr)\langle[\B^N]\rangle
        & (p,q)=(a,b),\\
   0 & \text{otherwise}.
 \end{cases}
\end{equation}
\end{fact}
\end{minipage}\par

\begin{proof}[Proof sketch]
Use the standard double-complex spectral sequence \cite[Lemmas~12.25.1 and 12.25.3]{Stacks} and the cohomology of a Koszul complex for a regular sequence \cite[Corollary~4.8.4(i)]{Schapira}. By definition, \(E_0^{p,q}=K_R^{p,q}\). Applying Fact~\ref{fact:koszul-h} in the vertical direction gives \(E_1^{p,q}=H^q(K_R^{p,*},d_v)\) and hence \eqref{eq:koszul-e1}. Applying the fact to the horizontal differential induced on \(E_1\) gives \eqref{eq:koszul-e2}.
\end{proof}

\begin{remark}\label{rem:coefficient-change}
There is an isomorphism of \(\Q\)-linear cochain complexes

\begin{equation}\label{eq:koszul-de-rham}
 K_R^\bullet(t_1-1,\ldots,t_N-1)\otimes_R M
 \cong\Omega^\bullet_{A/\Q}.
\end{equation}

In the commuting-operator description of \cite[Section~4.8]{Schapira}, take \(M=A\) and \(n=N\). Give \(M\) its \(R\)-module structure by

\begin{equation*}
 t_i\cdot f=(1+\partial_i)f,\qquad
 \partial_i=\frac{\partial}{\partial x_i}.
\end{equation*}

The operators are then given by

\begin{equation*}
 \varphi_i(f)=(t_i-1)\cdot f=\partial_i f .
\end{equation*}

This gives \eqref{eq:koszul-de-rham}. Compare \cite[Example~4.8.5(ii), p.~77]{Schapira}. With this choice of coefficients, the spectral sequences for the Koszul and de Rham complexes with the same filtration are isomorphic on each page.
\end{remark}

\section*{Acknowledgments}
This work was supported by JSPS KAKENHI Grant Number JP25K06999.

\section*{AI disclosure}
ChatGPT (OpenAI) was used to assist with translation of the author's Japanese notes into English and English-language editing, redrawing the author's hand-drawn originals in TikZ (Figures~\ref{fig:wb}, \ref{fig:example-flow}, and~\ref{fig:curve-orientations}), routine typesetting, and checks of mathematical notation and internal consistency in preparing this work. It also provided assistance in locating references and clarifying classical results. All mathematical statements, proofs, references, and the final text were independently checked by the author, who takes full responsibility for the work.

\end{document}